\documentclass[12pt,reqno]{amsart}

\usepackage{amsfonts,color,amsmath,amssymb}

\def\Box{\vcenter{\vbox{\hrule\hbox{\vrule
     \vbox to 8.8pt{\hbox to 10pt{}\vfill}\vrule}\hrule}}}

\newcommand{\tr}{\textup{Tr}}

\def\sgn{\textup{sgn}}
\newcommand{\ra}{\rangle}
\newcommand{\la}{\langle}
\newcommand{\F}{{\mathbb F}}

\newcommand{\cQ}{{\mathcal Q}}

\newcommand{\cP}{{\mathcal P}}
\newcommand{\cM}{{\mathcal M}}
\newcommand{\cL}{{\mathcal L}}
\newcommand{\PG}{\textup{PG}}
\newcommand{\PGO}{\textup{PGO}}

\newtheorem{thm}{Theorem}

\newtheorem{lemma}[thm]{Lemma}

\numberwithin{equation}{section}

\newtheorem{remark}[thm]{Remark}

\begin{document}
\title[$m$-ovoids of $Q(4,q)$]{A Family of $m$-ovoids of parabolic quadrics}

\author[Feng, Momihara and Xiang]{Tao Feng$^*$, Koji Momihara$^{\dagger}$ and Qing Xiang}
\thanks{$^*$Research supported in part by Fundamental Research Fund for the Central Universities of China, the National Natural Science Foundation of China under Grant No. 11201418 and 11422112.}
\thanks{$^{\dagger}$
Research supported by JSPS under Grant-in-Aid for Young Scientists (B) 25800093 and Scientific Research (C) 24540013.}

\address{Department of Mathematics, Zhejiang University, Hangzhou 310027, Zhejiang, P. R. China}
\email{tfeng@zju.edu.cn}

\address{Faculty of Education, Kumamoto University, 2-40-1 Kurokami, Kumamoto 860-8555, Japan} \email{momihara@educ.kumamoto-u.ac.jp}

\address{Department of Mathematical Sciences, University of Delaware, Newark, DE 19716, USA} \email{xiang@math.udel.edu}

\begin{abstract}
We construct a family of $\frac{(q-1)}{2}$-ovoids of $Q(4,q)$, the parabolic quadric of $\PG(4,q)$, for $q\equiv 3\pmod 4$. The existence of $\frac{(q-1)}{2}$-ovoids of $Q(4,q)$ was only known for $q=3, 7,$ or $11$. Our construction provides the first infinite family of $\frac{(q-1)}{2}$-ovoids of $Q(4,q)$. Along the way, we also give a construction of $\frac{(q+1)}{2}$-ovoids in $Q(4,q)$ for $q\equiv 1\pmod{4}$.
\end{abstract}

\maketitle

\section{Introduction}
This is a paper about $m$-ovoids of classical polar spaces. More specifically we are concerned with $m$-ovoids of $Q(4,q)$, the parabolic quadric of $\PG(4,q)$. Tight sets and $m$-ovoids are important substructures of classical polar spaces.  They are not only interesting in their own right, but also can give rise to other geometric/combinatorial objects, such as translation planes, strongly regular graphs, two-weight codes. The concept of a tight set was first defined by Payne \cite{payne} for generalized quadrangles, and it was extended to finite classical polar spaces by Drudge \cite{drudge}. The notion of an $m$-ovoid came from ``ovoids" of the projective space $\PG(3,q)$. It was Thas \cite{thaslaa} who first defined $m$-ovoids for generalized quadrangles, and later Shult and Thas \cite{ST1} extended the notion of an $m$-ovoid to finite classical polar spaces. In \cite{bambergCombin}, the authors unified the notion of a tight set and that of an $m$-ovoid of generalized quadrangles by defining intriguing sets of generalized quadrangles. In a subsequent paper \cite{bambergjcta}, the authors extended the concept of an intriguing set to finite classical polar spaces.

A finite {\it generalized quadrangle} (GQ) of order $(s,t)$ is an incidence structure $\cQ=(\cP, \cL)$ of points and lines such that two distinct lines intersect in at most one point, every line is incident with exactly $s+1$ points and every point is incident with exactly $t+1$ lines, and most importantly, that given any non-incident point-line pair $(P, L)$, $P\in\cP$ and $L\in\cL$, there is a unique line incident with $P$ and concurrent with $L$. Let $m$ be a positive integer. A subset $\cM$ of points of a generalized quadrangle $\cQ$ is called an {\it $m$-ovoid} if every line of $\cQ$ meets $\cM$ in exactly $m$ points. (Without giving the details we mention that $m$-ovoids of classical polar spaces can be defined similarly.) A 1-ovoid will be simply called an {\it ovoid}. We will only consider classical generalized quadrangles, namely, $W(3,q)$, $Q(4,q)$, $H(3,q^2)$, $Q^-(5,q)$, and $H(4,q^2)$, which are classical polar spaces of rank 2. In particular, we will concentrate on $Q(4,q)$, the parabolic quadric in $\PG(4,q)$, which is a generalized quadrangle of order $(q,q)$; the points of $Q(4,q)$ are the points of a non-singular quadric in $\PG(4,q)$, and the lines of $Q(4,q)$ are the lines contained in that quadric. Ovoids of $Q(4,q)$ are of special importance since from them we obtain spreads of $W(3,q)$, which in turn give rise to symplectic translation planes.

Ovoids of finite classical polar spaces are rare. They tend to exist only in classical polar spaces of low rank. In \cite{bambergCombin}, a systematic treatment of $m$-ovoids and tight sets in finite generalized quadrangles was given. It was commented in the same paper that few constructions of $m$-ovoids are known. In Example 5 of \cite{bambergCombin}, the authors gave a number of examples of $m$-ovoids of $Q(4,q)$, with $q=5,7,9,$ or $11$, found by a computer. In this paper, we generalize the examples in Example 5 of \cite{bambergCombin} into an infinite family in the case where $q\equiv 3\pmod 4$. As an illustration of our method, we also give a construction of a family of $\frac{q+1}{2}$-ovoids of $Q(4,q)$ when $q\equiv 1\pmod 4$. Note that there exist $\frac{q+1}{2}$-ovoids of $Q(4,q)$ for all odd $q$ (see Section 3 of \cite{bambergCombin}). So the $\frac{q+1}{2}$-ovoids obtained in this paper are of less importance. Nevertheless, we include a short section on $\frac{q+1}{2}$-ovoids of $Q(4,q)$ since the construction is direct and is done using the same idea as in the construction of $\frac{q-1}{2}$-ovoids of $Q(4,q)$. Our approach to the construction of these $m$-ovoids is similar to the one used previously in \cite{fmx}, \cite{jande}, and \cite{morgan}: we prescribe an automorphism group for the $m$-ovoids that we intend to construct, and then take unions of orbits of the point set of $Q(4,q)$ under the action of the prescribed automorphism group. Of course this approach had been used previously in the constructions of $m$-ovoids, tight sets, and other geometric objects. In most of the known constructions, e.g. \cite{CP}, usually the prescribed automorphism group is so large that there are only a limited number of orbits under the action of the group; as a consequence, it is possible to study the parameters of the related tactical decomposition explicitly. 

However, it is unusual for the geometric objects that we intend to construct to have very large automorphism groups. Erd\H{o}s and R\'enyi \cite{cameron} had shown that almost all graphs have no nontrivial automorphism, and  this is even true for various special classes of graphs \cite{er, wormald}. While we do not claim that almost all of the geometric objects that we are interested in have trivial automorphism groups, we believe that most of them should have fairly small automorphism groups.  So the new ingredient in our approach to the construction of $m$-ovoids is that we prescribe an automorphism group of medium size.  Consequently the number of orbits of the action of the group is large, and geometric argument for analyzing the intersection of the orbits with hyperplanes seems impossible. In our case, we succeed in deriving an algebraic expression of the underlying set, and then give a purely algebraic proof. The difficult part lies in the derivation of the algebraic expression which is valid for all prime power $q$ congruent to 3 modulo 4. It requires a good understanding of the geometric structure and a translation into algebraic language. We accomplish this by examining the interplay of the putative $m$-ovoids with some special hyperplanes.

\section{A family of $\frac{(q-1)}{2}$-ovoids of $Q(4,q)$}

In this section, we generalize some examples of $m$-ovoids of $Q(4,q)$ found by a computer in \cite{bambergCombin} into an infinite family. As we already mentioned in the introduction,  the basic idea is to prescribe an automorphism group of the $m$-ovoids to be constructed. The choice of the automorphism group is somewhat evident from the examples compiled in Example 5 of Section 8 in \cite{bambergCombin}. The difficulty lies in choosing orbits of the point set of $Q(4,q)$ under the action of the prescribed group so that the union of the chosen orbits forms an $m$-ovoid of $Q(4,q)$.

\subsection{The model of $Q(4,q)$ and the group $A$}
Let $q\equiv 3\pmod{4}$ be a prime power, and let $V=\F_q\times\F_q\times\F_{q^2}\times \F_q$. We view $V$ as a 5-dimensional vector space over $\F_q$, and for any $(x_1,x_2,y,z)\in V$, define
$$f(x_1,x_2,y,z)=x_1x_2+y^{q+1}-z^2.$$
It is easy to see that $f$ is a non-degenerate quadratic form on $V$, hence the set of zeros of $f$,  $\{(x_1,x_2,y,z)\in V\mid  x_1x_2+y^{q+1}-z^2=0\}$, gives a parabolic quadric of $\PG(4,q)$. This is the model that we will be using for $Q(4,q)$. Formally, we define 
$$Q(4,q)=\{\langle (x_1,x_2,y,z)\rangle\mid 0\neq (x_1,x_2,y,z)\in V, \; f(x_1,x_2,y,z)=0\}.$$
In the rest of the paper, in order to simplify notation, we often simply write $(x_1,x_2,y,z)$ for the projective point $\langle(x_1,x_2,y,z)\rangle$ of $\PG(4,q)$. For future use, we note that the polar form $B$ of $f$ is defined by
$$B((x_1,x_2,y,z), (x'_1,x'_2,y',z'))=x'_1x_2+x_1x'_2+y^qy'+y'^qy-2zz', \;\forall (x_1,x_2,y,z),  (x'_1,x'_2,y',z')\in V.$$

Let $\square_q$ and $\blacksquare_q$ be the set of nonzero squares of $\F_q$ and the set of nonsquares of $\F_q$, respectively. We define $H$ to be the following set of mappings from $\PG(4,q)$ to itself.
\[H=\{
T_{\beta,w}: (x_1,x_2,y,z)\mapsto (x_1\beta  ,x_2\beta^{-1} ,wy,w^{\frac{q+1}{2}}z )\mid \beta\in\square_{q},\,w^{q+1}=1, w\in\F_{q^2}\}.
\]

\begin{lemma}
The set $H$ defined above forms a cyclic subgroup of order $\frac{q^2-1}{2}$ of $\PGO(5,q)$.
\end{lemma}
\proof 
It is clear that $T_{\beta,w}$ is $\F_q$-linear since the image of $(x_1,x_2,y,z)$ under $T_{\beta,w}$ is obtained by 
multiplying each coordinate of $(x_1,x_2,y,z)$ by some element of $\F_{q^2}$. Furthermore, $T_{\beta,w}$ preserves the quadratic form $f$  since 
\[
f(T_{\beta,w}(x_1,x_2,y,z))=f(x_1\beta,x_2\beta^{-1},wy,w^{\frac{q+1}{2}}z)=x_1x_2+y^{q+1}-z^2=f(x_1,x_2,y,z)
\] 
by the assumption that $w^{q+1}=1$. 
Note that since $q\equiv 3\pmod 4$, we have $\gcd(q+1, \frac{q-1}{2})=1$. Hence, the set $H$ is actually a cyclic subgroup of order $\frac{q^2-1}{2}$  of $\PGO(5,q)$.
\qed

\vspace{0.3cm}
We define two further involutions of $\PGO(5,q)$ as follows:
\begin{align*}
   &\sigma:\,(x_1,x_2,y,z)\mapsto (x_1,x_2,-y^q,z),\\
   &\tau:\,(x_1,x_2,y,z)\mapsto(x_2,x_1,y,-z).
\end{align*}
It is straightforward to check that $\sigma\tau=\tau\sigma$. We now define $A$ to be the subgroup of $\PGO(5,q)$ generated by $H, \sigma$, and $\tau$. That is, $A:=\la H,\sigma,\tau\ra$. The group $A$ has order $2(q^2-1)$, and it is isomorphic to $C_{\frac{q^2-1}{2}}\rtimes(C_2\times C_2)$, where $C_m$  denotes a cyclic group of order $m$ for any positive integer $m$.

The group $A$ acts on the set of points of $Q(4,q)$. Below we collect some basic observations about this action.  Given a point $P$ of $\PG(4,q)$, we will use $O(P)$ to denote the orbit in which $P$ lies under the action of $A$.
\begin{enumerate}
\item The plane $y=0$ intersects $Q(4,q)$ in a conic $C=\{\langle (x_1,x_2,0,z)\rangle \mid 0\neq (x_1,x_2,0,z)\in V, \; x_1x_2=z^2\}$. The group $A$ interchanges the two points $e_1=(1,0,0,0)$ and $e_2=(0,1,0,0)$, and acts transitively on the remaining $(q-1)$ points of $C$.
\item The hyperplane $z=0$ intersects $Q(4,q)$ in an elliptic quadric $Q^-(3,q)=\{\langle (x_1,x_2,y,0)\rangle \mid  0\neq (x_1,x_2,y,0)\in V,\; x_1x_2+y^{q+1}=0\}$. From the definition of $A$, one can see that $A$ fixes $Q^-(3,q)$ setwise, and $Q^-(3,q)$ is partitioned into three orbits under the action of $A$. The three orbits are $\{e_1, e_2\}$, $O(1,-1,1,0)$ and $O(1,1,\mu,0)$, where $\mu\in\F_{q^2}$ and $\mu^{q+1}=-1$.
\end{enumerate}

We now give a summary of the orbit structure of $Q(4,q)$ under the action of $A$. There is a unique orbit $\{e_1,e_2\}$ of length 2, a unique orbit $O(1,1,0,1)=C\setminus \{e_1,e_2\}$ of length $q-1$, and a unique orbit of length $q+1$, namely, $O(0,0,1,1)$. All other orbits have length either $\frac{q^2-1}{2}$ or $q^2-1$. There are $\frac{q+1}{2}$ orbits of length $\frac{q^2-1}{2}$, which are given below.
\begin{enumerate}
\item The set of points $(x_1,x_2,y,0)$ of $Q(4,q)$ with $y\ne 0$ splits into two orbits,  $O(1,1,\mu,0)$ and $O(1,-1,1,0)$, where $\mu$ is the same as above. Both orbits have length $\frac{q^2-1}{2}$.
\item $O(1,-t^2,y,1)$, $1+t^2\in\square_q$, $y\in \F_{q^2}^*$, and $y^{q+1}=1+t^2$. There are a total of 
$\frac{q-3}{2}$ such orbits, each having length $\frac{q^2-1}{2}$.
\end{enumerate}

There are a total of $\frac{3(q+1)}{4}$ orbits of length $q^2-1$, which we describe now.
\begin{enumerate}
\item $O(1,0,1,1)$ and $O(-1,0,1,1)$,
\item $O(x_1,x_2,y,1)$ with $x_1x_2\neq 0$ and $y\neq 0$, and either $x_1x_2\in\square_q$ or $1-x_1x_2\in\blacksquare_q$.
\end{enumerate}

In the rest of this section, an orbit of length $\frac{q^2-1}{2}$ will be referred to as a short orbit, and an orbit of length $q^2-1$ as a long orbit. Another way to organize the orbits of the points of $Q(4,q)$ under the action of $A$ is as follows. There are 7 orbits in which a representative has a coordinate being zero: these are $O(1,0,0,0)=\{e_1,e_2\}$, $O(1,1,0,1)=C\setminus\{e_1,e_2\}$, $O(0,0,1,1)$, $O(1,1,\mu,0)$, $O(1,-1,1,0)$, $O(1,0,1,1)$, and $O(-1,0,1,1)$. All the other orbits have the form $O(x_1,x_2,y,1)$ with $x_1x_2y\neq 0$. Among these orbits, $O(x_1,x_2,y,1)$ is a short orbit if and only if $x_1x_2\in\blacksquare_q$ and $1-x_1x_2\in\square_q$. (In the above we parameterized the short orbits $O(x_1,x_2,y,1)$ more explicitly so that we can count them efficiently.) We state this claim as a lemma.

\begin{lemma}
Let $x_1,x_2\in\F_q^*$, $y\in\F_{q^2}^*$, and $(x_1,x_2,y,1)\in Q(4,q)$. Then the $A$-orbit $O(x_1,x_2,y,1)$ is a short orbit if and only if $x_1x_2\in\blacksquare_q$ and $1-x_1x_2\in\square_q$. 
\end{lemma}

\proof We will compute the point stabilizer ${\rm Stab}_{A}(x_1,x_2,y,1)$ of $(x_1,x_2,y,1)$ in $A=\la H,\sigma,\tau\ra$. First, if $T_{\beta,w}(x_1,x_2,y,1)=(x_1,x_2,y,1)$, then there exists $\theta\in\F_q^*$ such that $(\beta x_1, \beta^{-1}x_2, wy, w^{\frac{q+1}{2}})=\theta(x_1,x_2,y,1)$. It follows that $\beta=\theta=\beta^{-1}$, $w=\theta=w^{\frac{q+1}{2}}$. Since $\beta\in\square_q$, we have $\theta=1$. Hence $\beta=w=1$. So the conclusion is that if $T_{\beta,w}(x_1,x_2,y,1)=(x_1,x_2,y,1)$, then $T_{\beta,w}=T_{1,1}$. Secondly, if $T_{\beta,w}\cdot \sigma(x_1,x_2,y,1)=(x_1,x_2,y,1)$, then there exists $\theta\in\F_q^*$ such that $(\beta x_1, \beta^{-1}x_2, -wy^q, w^{\frac{q+1}{2}})=\theta(x_1,x_2,y,1)$. It follows that $\beta=\theta=\beta^{-1}$, $-wy^q=y$, and $\theta=w^{\frac{q+1}{2}}$. Again we must have $\theta=1$ since $\beta\in\square_q$. Also since $(x_1,x_2,y,1)\in Q(4,q)$, we have $1-x_1x_2=y^{q+1}$. And $(1-x_1x_2)^{\frac{q-1}{2}}=y^{(q+1)\cdot \frac{q-1}{2}}=(-\frac{1}{w})^{\frac{q+1}{2}}=1$. It follows that $1-x_1x_2\in\square_q$. So the conclusion is that $T_{\beta,w}\cdot \sigma(x_1,x_2,y,1)=(x_1,x_2,y,1)$ if and only if $\beta=1$, $1-x_1x_2\in\square_q$, and $w=-y^{1-q}$. 

Similarly, we have that $T_{\beta,w}\cdot \tau(x_1,x_2,y,1)=(x_1,x_2,y,1)$ if and only if $w=-1$, $x_1x_2\in\blacksquare_q$, and $\beta=-\frac{x_1}{x_2}$. 
We now see that $|{\rm Stab}_{A}(x_1,x_2,y,1)|=4$ if and only if $x_1x_2\in\blacksquare_q$ and $1-x_1x_2\in\square_q$. The proof of the lemma is complete.
\qed

\subsection{The description of the set $\cM$}
We are now ready to give the promised construction of $\frac{(q-1)}{2}$-ovoids of $Q(4,q)$. As we mentioned before, the $\frac{(q-1)}{2}$-ovoids that we intend to construct are unions of orbits of the action of $A$ on $Q(4,q)$. Let $a$ be a fixed element of $\F_q^*$ such that $1+a^2\in\square_q$, $d\in\F_q^*$ be a fixed element such that $d^2=a^{-2}+1$, $\mu$ be the same as above (that is, $\mu\in\F_{q^2}$ and $\mu^{q+1}=-1$), and let $S=\{\langle (x_1,x_2,y,1)\rangle \in Q(4,q)\mid 1+a^2x_1x_2 \in\square_q, x_1x_2y\neq 0\}$.  In a compact way, we define 
\begin{equation}\label{movoid1}
\cM=S\cup O(1,0,1,1)\cup O(1,1,\mu,0)\cup O(1,1,0,1)\cup O(-1,1,ad,a).
\end{equation}

We remark that the set $S$ is invariant under the action of $A$. So $S$ is a union of $A$-orbits. The size of $S$ can be computed easily. For $(x_1,x_2,y,1)\in S$, we must have $x_1x_2\in \square_q-a^{-2}$ and $x_1x_2\neq 0$ or $1$ (note that if $x_1x_2=1$, since $x_1x_2+y^{q+1}-1=0$, we have $y=0$, which is impossible). So there are $\frac{q-1}{2}-2=\frac{q-5}{2}$ choices for $x_1x_2$; for each such choice of $x_1x_2$, since  $y^{q+1}=1-x_1x_2$, we have $q+1$ choices for $y$. Therefore $|S|=\frac{(q-5)(q^2-1)}{2}.$  It follows that 
$$|\cM|=\frac{(q-5)}{2}(q^2-1)+(q^2-1)+2\cdot \frac{q^2-1}{2}+(q-1)=\frac{(q-1)}{2}(q^2+1).$$

The geometric interpretation of $\cM$ is as follows (this is actually how we chose the orbits to form $\cM$ {\color{red}in} the first place.) Let $R=(a,-a,0,1)$. Then $R^{\perp}=\{(x_1,x_2,y,z)\mid z=\frac{a}{2}(x_2-x_1)\},$ where $\perp$ is defined by the bilinear form $B$, i.e., the polar form of $f$. We have
\begin{enumerate}
\item $\cM$ contains three orbits in which a representative has a zero coordinate: \[O(1,0,1,1) \textup{--long},\; O(1,1,\mu,0) \textup{--short},\; O(1,1,0,1).\]
\item $\cM$ contains all the long orbits that have nonempty intersection with $R^\perp$;
\item $\cM$ contains all the short orbits $O(b,-b,c,1)$ such that $R^\perp\cap O(\pm b,\mp b,c,1)\neq \emptyset$ and the intersection contains a ``rational" point (here, a point $(x_1,x_2,y,z)$ is called rational if $y\in\F_q$).
\end{enumerate}

We are now ready to state the main theorem of this section.
\begin{thm}\label{main1}
The above set $\cM$ is a $\frac{(q-1)}{2}$-ovoid of $Q(4,q)$.
\end{thm}

\proof We will show that every line of $Q(4,q)$ meets $\cM$ in exactly $\frac{q-1}{2}$ points. Any line of $Q(4,q)$ takes the form $\overline{PQ}=\{\lambda P+Q \mid P,Q\in Q(4,q), P\perp Q, \lambda\in\F_q\cup \{\infty\}\}$, where $\infty\cdot P+Q=P$; in order to simplify the computations, we will take $P=(x_1,x_2,y,z)\in Q(4,q)$ with $x_2=0$. Note that this is always possible since we can vary the point $P$ on the line $\overline{PQ}$. Also since $\cM$ is fixed setwise by $A$, we see that $|\overline{g(P)g(Q)}\cap\cM|=|\overline{PQ}\cap \cM|$ for any line $\overline{PQ}$ of $Q(4,q)$ and any $g\in A$. As we mentioned above, there are 7 orbits of the points of $Q(4,q)$ in which a representative has a coordinate being zero. By examining the 7 orbits,we see that it suffices to consider the following four choices for $P$, namely, $P=(1,0,0,0)$, $(1,0,1,1)$, $(-1,0,1,1)$, or $(0,0,1,1)$.

\begin{enumerate}
\item $P=(1,0,0,0)$. Choose $Q=(a_1,a_2,b,c)\in Q(4,q)$ such that $P\perp Q$. We have $a_2=0$ and $b^{q+1}=c^2$, where $b\in \F_{q^2}$ and $c\in \F_q^*$. We will compute $|\overline{PQ}\cap \cM|$. Note that $\overline{PQ}=\{(a_1+\lambda,0,b,c)\mid \lambda\in \F_q\cup \{\infty\}\}$. First we claim that $\overline{PQ}\cap S=\emptyset$. This is clear since any point in $S$ has the second coordinate nonzero. For the same reason, we have that $\overline{PQ}\cap O(1,1,\mu, 0)=\emptyset$, $\overline{PQ}\cap O(1,1,0,1)=\emptyset$, and $\overline{PQ}\cap O(-1,1,ad,a)=\emptyset$. Next, we claim that $|\overline{PQ}\cap O(1,0,1,1)|=\frac{q-1}{2}$. To prove the claim, we first note that the orbit $O(1,0,1,1)$ has length $q^2-1$, and it consists of the following points, $(\beta,0,w,w^{\frac{q+1}{2}})$ and $(0,\beta,w,-w^{\frac{q+1}{2}})$, where $\beta\in\square_q$ and $w^{q+1}=1$, $w\in\F_{q^2}$. It is then easy to see that $\overline{PQ}\cap O(1,0,1,1)=\overline{PQ}\cap \{(\beta,0,w,w^{\frac{q+1}{2}})\mid \beta\in\square_q, w^{q+1}=1, w\in\F_{q^2}\}$. We will find that the latter intersection has size $\frac{q-1}{2}$. To see this, set $(a_1+\lambda,0,b,c)=\theta_{\lambda} (\beta,0,w,w^{\frac{q+1}{2}})$, where $\theta_{\lambda}\in\F_q^*$. We have 
$$a_1+\lambda=\theta_{\lambda}\beta,\; b=\theta_{\lambda}w,\; c=\theta_{\lambda}w^{\frac{q+1}{2}}.$$
It follows that $w^{\frac{q-1}{2}}=c/b$. Since $\gcd(q+1,\frac{q-1}{2})=1$ and $w^{q+1}=1$, we have that $w=(c/b)^s$ for some integer $s$, and $\theta_{\lambda}=w^{-\frac{q+1}{2}}c$. This shows that as $\lambda$ runs through the elements of $\F_q$, $\theta_{\lambda}$ is independent of $\lambda$. Now it is clear that $a_1+\lambda=\theta_{\lambda}\beta$ has $\frac{q-1}{2}$ solutions in $\lambda$ since $\beta\in\square_q$. Therefore we have  shown that $|\overline{PQ}\cap \{(\beta,0,w,w^{\frac{q+1}{2}})\mid \beta\in\square_q, w^{q+1}=1, w\in\F_{q^2}\}|=\frac{q-1}{2}$, from which it follows that $|\overline{PQ}\cap O(1,0,1,1)|=\frac{q-1}{2}.$ Now $|\overline{PQ}\cap {\mathcal M}|=|\overline{PQ}\cap O(1,0,1,1)|=\frac{q-1}{2}.$

\vspace{0.1in}

\item $P=(1,0,1,1)$. Choose $Q=(a_1,a_2,b,c)\in Q(4,q)$ such that $P\perp Q$. By varying $Q$ on the line $\overline{PQ}$, we may assume that $a_1=0$. Since $Q\in Q(4,q)$, we have $b^{q+1}=c^2$, where $b\in \F_{q^2}$ and $c\in \F_q$. Also $a_2+b^q+b-2c=0$ since $P\perp Q$. We will assume that $a_2\neq 0$ and $b\not\in \F_q$ (the case where $a_2=0$ or $b\in\F_q$ is easy to handle). Consequently $c\neq 0$, and we may further assume that $c=1$. So the conditions satisfied by $a_2,b$ are
$$b^{q+1}=1, b\in\F_{q^2}\setminus\F_q,\; a_2+b^q+b-2=0, a_2\in\F_q^*.$$
We will compute $|\overline{PQ}\cap \cM|$, where $\overline{PQ}=\{(\lambda,a_2,b+\lambda,1+\lambda)\mid \lambda\in \F_q\cup \{\infty\}\}$. We claim that 
\begin{equation}\label{empty}
\overline{PQ}\cap O(1,1,0,1)=\emptyset,\end{equation}
\begin{equation}\label{onepoint}
\overline{PQ}\cap O(1,0,1,1)=\{P\},\end{equation}
\begin{equation}\label{zeroorone}
|\overline{PQ}\cap O(1,1,\mu,0)|=[[-a_2\in\square_q]],\end{equation}
\begin{equation}\label{large}|\overline{PQ}\cap S|=\frac{q-3}{2}-\frac{|\{\lambda\in\F_q\mid \lambda^2+(2+a^2a_2)\lambda+1=0\}|}{2}-[[-a_2\in\square_q]],\end{equation}
\begin{equation}\label{0or1}|\overline{PQ}\cap O(-1,1,ad,a)|=\frac{|\{\lambda\in\F_q\mid \lambda^2+(2+a^2a_2)\lambda+1=0\}|}{2},\end{equation}
where $[[-a_2\in\square_q]]$ is the Kronecker delta function taking value $1$ if $-a_2\in\square_q$ and value $0$ otherwise. Once we prove all five claims above, we will see that $|\overline{PQ}\cap \cM|=\frac{q-1}{2}.$ First we consider $\overline{PQ}\cap O(1,1,0,1)$. Recall that $O(1,1,0,1)=C\setminus\{e_1,e_2\}$ is an orbit of length $q-1$, where $C$ is the conic $\{\langle (x_1,x_2,0,z)\rangle \mid 0\neq (x_1,x_2,0,z)\in V, \; x_1x_2=z^2\}$.  If $(x_1,x_2,0,z)=\theta_{\lambda}(\lambda,a_2,b+\lambda,1+\lambda)$ for some $\theta_{\lambda}\in\F_q^*$, then $b=-\lambda\in\F_q$. But we already assumed that $b\not\in \F_q$. This proves that $\overline{PQ}\cap O(1,1,0,1)=\emptyset$. 

Next, we consider $\overline{PQ}\cap O(1,0,1,1)$. It is clear that $P\in \overline{PQ}\cap O(1,0,1,1)$. So we only need to consider $\{(\lambda,a_2,b+\lambda,1+\lambda)\mid \lambda\in \F_q\}\cap O(1,0,1,1)$. Note that $O(1,0,1,1)$ has length $q^2-1$, and it consists of the following points, $(\beta,0,w,w^{\frac{q+1}{2}})$ and $(0,\beta,w,-w^{\frac{q+1}{2}})$, where $\beta\in\square_q$ and $w^{q+1}=1$, $w\in\F_{q^2}$. Since $a_2\neq 0$, we see that $\{(\lambda,a_2,b+\lambda,1+\lambda)\mid \lambda\in \F_q\}\cap \{(\beta,0,w,w^{\frac{q+1}{2}})\mid \beta\in\square_q, w^{q+1}=1\}=\emptyset.$ We claim that also $\{(\lambda,a_2,b+\lambda,1+\lambda)\mid \lambda\in \F_q\}\cap \{(0,\beta,w,-w^{\frac{q+1}{2}})\mid \beta\in\square_q, w^{q+1}=1\}=\emptyset.$ To see this, let $\theta_{\lambda}(\lambda,a_2,b+\lambda,1+\lambda)=(0,\beta,w,-w^{\frac{q+1}{2}})$ for some $\theta_{\lambda}\in\F_q^*$. We have $\lambda=0$, and the only possible point in the intersection $\{(\lambda,a_2,b+\lambda,1+\lambda)\mid \lambda\in \F_q\}\cap \{(0,\beta,w,-w^{\frac{q+1}{2}})\mid \beta\in\square_q, w^{q+1}=1\}$ is the point $Q$. Now $Q=(0,a_2,b,1)$ is in the aforementioned intersection if and only if there exists some $\theta_0\in\F_q$ such that $\theta_0 a_2=\beta, \theta_0 b=w, \theta_0=-w^{\frac{q+1}{2}}\in \{1,-1\}$, which, in turn, imply that $a_2 b^{\frac{q+1}{2}}\in \blacksquare_q$. We claim that in fact $a_2 b^{\frac{q+1}{2}}\in \square_q$, from which it follows that $Q$ is not in the intersection and therefore $\{(\lambda,a_2,b+\lambda,1+\lambda)\mid \lambda\in \F_q\}\cap \{(0,\beta,w,-w^{\frac{q+1}{2}})\mid \beta\in\square_q, w^{q+1}=1\}=\emptyset.$ For $x\in \F_q$, we define the sign of $x$, $\sgn(x)\in \F_q$, by
\begin{equation}\label{eqn:sgn}
\sgn(x)
=\begin{cases} 1,\quad &\textup{if $x\in\square_q$},\\
-1,\quad & \textup{if $x\in\blacksquare_q$},\\
0,\quad & \textup{if $x=0$}.\end{cases}
\end{equation}
Since $b^{q+1}=1$ and $b+b^q=2-a_2$, we have $4-a_2=2+b+b^q=(1+b)(1+b^q)=b^{-1}(1+b)^2$. It follows that $(4-a_2)^{\frac{q-1}{2}}=b^{\frac{q+1}{2}}$, and $\sgn(4-a_2)=\sgn(b^{\frac{q+1}{2}})$. On the other hand, since $b$ and $b^q$ are two roots of $X^2-(2-a_2)X+1$ in $\F_{q^2}$ and $b\not\in\F_q$, we have $(2-a_2)^2-4=a_2^2-4a_2\in\blacksquare_q$, and it follows that $\sgn(a_2)\sgn(4-a_2)=1$. Hence $\sgn(a_2)\sgn(b^{\frac{q+1}{2}})=1$. That is, $a_2b^{\frac{q+1}{2}}\in \square_q$. We have finished the proof of \eqref{onepoint}.

We now consider $\overline{PQ}\cap O(1,1,\mu,0)$. It is clear that the only possible point in the intersection is $(-1,a_2,b-1,0)$, and $(-1,a_2,b-1,0)\in \overline{PQ}\cap O(1,1,\mu,0)$ if and only if $-a_2\in\square_q$. This proves \eqref{zeroorone}.

We consider $\overline{PQ}\cap S$ next. Recall that $S=\{\langle (x_1,x_2,y,1)\rangle \in Q(4,q)\mid 1+a^2x_1x_2 \in\square_q, x_1x_2y\neq 0\}$.  Take a point  $(x_1,x_2,y,1)\in S$ and set $(x_1,x_2,y,1)=\theta_{\lambda}(\lambda,a_2,b+\lambda,1+\lambda)$ for some $\theta_{\lambda}\in\F_q^*$. We have 
\begin{align*}
x_1=&\,\theta_{\lambda} \lambda,\\
x_2=&\,\theta_{\lambda} a_2,\\
y=&\,\theta_{\lambda}(b+\lambda),\\
1=&\,\theta_{\lambda} (1+\lambda).
\end{align*}
It follows that $1+a^2x_1x_2=1+\theta_{\lambda}^2a^2\lambda a_2=1+(1+\lambda)^{-2}a^2\lambda a_2$. Since $(x_1,x_2,y,1)\in S$, we must have $(1+\lambda)^2+a^2\lambda a_2\in \square_q$. Conversely, for each $\lambda\in\F_q$ satisfying $(1+\lambda)^2+a^2\lambda a_2\in \square_q$, as long as $\lambda\neq 0$ or $-1$, we can always solve the above system so that the point $(\lambda,a_2,b+\lambda,1+\lambda)$ belongs to $S$. Therefore, 
\begin{equation}\label{eq:quadeq1}
|\overline{PQ}\cap S|=|\{\lambda\in \F_q\,|\,(1+\lambda)^2+a^2\lambda a_2\in \square_q, \lambda\not=0,-1\}|. 
\end{equation}
If $\lambda=0$, $(1+\lambda)^2+a^2\lambda a_2\in \square_q$, and if $\lambda=-1$, $(1+\lambda)^2+a^2\lambda a_2\in \square_q$ or $\blacksquare_q$ according as $-a_2\in \square_q$ or not. 
Hence, continuing from \eqref{eq:quadeq1}, we have 
\[
|\overline{PQ}\cap S|=\sum_{\lambda\in\F_q,\; (1+\lambda)^2+\lambda a^2a_2\neq 0}\frac{\eta((1+\lambda)^2+\lambda a^2a_2)+1}{2}-1-[[-a_2\in \square_q]]\]
\[\hspace{0.6in}=\sum_{\lambda\in\F_q}\frac{\eta((1+\lambda)^2+\lambda a^2a_2)+1}{2}-\frac{|\{\lambda\in\F_q\mid \lambda^2+(2+a^2a_2)\lambda+1=0\}|}{2}-1-[[-a_2\in \square_q]]\]
\[\hspace{0.6in}=\frac{q-2}{2}+\frac{1}{2}\sum_{\lambda\in\F_q}\eta(\lambda^2+(2+a^2a_2)\lambda+1)-\frac{|\{\lambda\in\F_q\mid \lambda^2+(2+a^2a_2)\lambda+1=0\}|}{2}-[[-a_2\in \square_q]],\]
where $\eta$ is the quadratic character of $\F_q$. The discriminant of the quadratic polynomial $X^2+(2+a^2a_2)X+1$ is $\Delta:=a^2a_2(4+a^2a_2)$. If $\Delta=0$, then $a_2=-\frac{4}{a^2}$, and $a_2(4-a_2)=-\frac{16(1+a^2)}{a^4}$. Since $b$ and $b^q$ are two roots of the quadratic polynomial $X^2-(2-a_2)X+1$ in $\F_{q^2}$ and $b\not\in\F_q$, we have $(2-a_2)^2-4=a_2^2-4a_2\in\blacksquare_q$, implying that $a_2(4-a_2)\in \square_q$.  Since $a_2(4-a_2)=-\frac{16(1+a^2)}{a^4}$, we therefore have $1+a^2\in\blacksquare_q$, contradicting the assumption $1+a^2\in \square_q$. We have shown that the discriminant of the quadratic polynomial $X^2+(2+a^2a_2)X+1$ is nonzero, and it follows that $\sum_{\lambda\in\F_q}\eta(\lambda^2+(2+a^2a_2)\lambda+1)=-1$ by \cite[Theorem~5.48]{LN97}. Consequently, $$|\overline{PQ}\cap S|=\frac{q-3}{2}-\frac{|\{\lambda\in\F_q\mid \lambda^2+(2+a^2a_2)\lambda+1=0\}|}{2}-[[-a_2\in\square_q]].$$ This proves \eqref{large}.

Finally we consider $\overline{PQ}\cap O(-1,1,ad,a)$. Recall that $O(-1,1,ad,a)$ is a short orbit, and it consists of points $(-\beta,\beta^{-1},adw,aw^{\frac{q+1}{2}})$, where $\beta\in\square_q$ and $w^{q+1}=1$. Set $\theta_{\lambda}(\lambda,a_2,b+\lambda,1+\lambda)=(-\beta,\beta^{-1},adw,aw^{\frac{q+1}{2}})$ for some $\theta_{\lambda}\in\F_q^*$. We have

\begin{equation}\label{beta}
-\beta=\theta_{\lambda} \lambda,\end{equation}
\begin{equation}\label{betainverse}\beta^{-1}=\theta_{\lambda} a_2,\end{equation}
\begin{equation}\label{third}adw=\theta_{\lambda}(b+\lambda),\end{equation}
\begin{equation}\label{last}aw^{\frac{q+1}{2}}=\theta_{\lambda} (1+\lambda).
\end{equation}

Multiplying \eqref{beta} and \eqref{betainverse}, we obtain $\theta_{\lambda}^2\lambda a_2=-1$. Square both sides of \eqref{last}, we have $a^2=\theta_{\lambda}^2(1+\lambda)^2$. Therefore we have 

\begin{equation}\label{quadcond}
(1+\lambda)^2+a^2a_2\lambda=0.
\end{equation}

Note that we have seen above that $\Delta\neq 0$. So a necessary condition for a point $(\lambda,a_2,b+\lambda,1+\lambda)\in \overline{PQ}$ to lie in the orbit $O(-1,1,ad,a)$ is that $\Delta:=a^2a_2(4+a^2a_2)\in\square_q$. It follows that if $\Delta\in\blacksquare_q$, then $\overline{PQ}\cap O(-1,1,ad,a)=\emptyset$. Now assume that $\Delta\in\square_q$. Then there are $\lambda_1\neq \lambda_2$ in $\F_q$ such that $(1+\lambda_i)^2+a^2a_2\lambda_i=0$ for $i=1,2$. We will show that exactly one of the two points $(\lambda_i,a_2,b+\lambda_i,1+\lambda_i)\in \overline{PQ}$, $i=1,2$, lies in $O(-1,1,ad,a)$. 

Since $\gcd(q+1, \frac{q-1}{2})=1$, we can find integers $s$ and $t$ such that $(q+1)t+\frac{(q-1)s}{2}=1$. Given $\lambda_i\in\F_q$ ($i=1$ or $2$) satisfying \eqref{quadcond}, we set
$$w_i=\left(\frac{d(1+\lambda_i)}{b+\lambda_i}\right)^s.$$
Since $b^{q+1}=1$, $b^q+b=2-a_2$, $d^2=1+a^{-2}$, and $\lambda_i$ satisfies \eqref{quadcond}, we have $(b+\lambda_i)^{q+1}=b^{q+1}+(b^q+b)\lambda_i+\lambda_i^2=1+(2-a_2)\lambda_i+\lambda_i^2=d^2(1+\lambda_i)^2$. It follows that the $w_i$'s defined above satisfy $w_i^{q+1}=1$, $w_i^{\frac{q-1}{2}}=\frac{d(1+\lambda_i)}{b+\lambda_i}$, and $\frac{w_i}{b+\lambda_i}\in \F_q$. Set $\theta_{\lambda_i}=\frac{adw_i}{b+\lambda_i}$, where $i=1, 2$. Then $\theta_{\lambda_i}\in\F_q$, and both \eqref{third} and \eqref{last} are satisfied. Next we show that exactly one of the following holds:  $\sgn(-\frac{adw_i}{b+\lambda_i})=1$, where $i=1, 2$. We achieve this by showing that 
$$\sgn(\frac{w_1}{b+\lambda_1}\cdot \frac{w_2}{b+\lambda_2})=-1.$$
It suffices to show that $(\frac{w_1}{b+\lambda_1}\cdot \frac{w_2}{b+\lambda_2})^{\frac{q-1}{2}}=-1$. Recall that $w_i^{\frac{q-1}{2}}=\frac{d(1+\lambda_i)} {b+\lambda_i}$, $i=1,2$, $\lambda_1\lambda_2=1$,  $\lambda_1+\lambda_2=-2-a_2a^2$, and $b^2+1=(2-a_2)b$. We see that 
$$\left(\frac{w_1}{b+\lambda_1}\cdot \frac{w_2}{b+\lambda_2}\right)^{\frac{q-1}{2}}=\frac{d^2(1+\lambda_1)
(1+\lambda_2)} {(b^2-(2+a^2a_2)b+1)^{\frac{q+1}{2}}}=\frac{-a_2(1+a^2)}{(-a_2(1+a^2)b)^{\frac{q+1}{2}}}$$
By assumption, $1+a^2\in \square_q$, we have $(1+a^2)^{\frac{q+1}{2}}=1+a^2$. Note that $\frac{q+1}{2}$ is even since $q\equiv 3\pmod 4$. We have
$$\left(\frac{w_1}{b+\lambda_1}\cdot \frac{w_2}{b+\lambda_2}\right)^{\frac{q-1}{2}}=\frac{-a_2}{(a_2b)^{\frac{q+1}{2}}}=\frac{-1}{a_2^{\frac{q-1}{2}}b^{\frac{q+1}{2}}}$$
Recall that in the proof of \eqref{onepoint}, we have shown that $a_2b^{\frac{q+1}{2}}\in\square_q$. (Also note that $b^{\frac{q+1}{2}}=\pm 1$.) Therefore $a_2^{\frac{q-1}{2}}b^{\frac{q+1}{2}}=1$. It follows that$\left(\frac{w_1}{b+\lambda_1}\cdot \frac{w_2}{b+\lambda_2}\right)^{\frac{q-1}{2}}=-1.$ We have now shown that exactly one of the following holds:  $\sgn(-\frac{adw_i}{b+\lambda_i})=1$, where $i=1, 2$. For the $i$ such that $\sgn(-\frac{adw_i}{b+\lambda_i})=1$, we define $\beta=-\frac{adw_i\lambda_i}{b+\lambda_i}$. Then $\beta\in\square_q$, and both \eqref{beta} and \eqref{betainverse} are satisfied. So in this case (i.e., the case where $\Delta\in\square_q$) we have $|\overline{PQ}\cap O(-1,1,ad,a)|=1$. In summary we have shown that $|\overline{PQ}\cap O(-1,1,ad,a)|=\frac{|\{\lambda\in\F_q\mid \lambda^2+(2+a^2a_2)\lambda+1=0\}|}{2}$.

\vspace{0.1in}
\item $P=(-1,0,1,1)$. This case is almost identical to the case where $P=(1,0,1,1)$. We omit the details.

\vspace{0.1in}
\item $P=(0,0,1,1)$. This case is rather easy. We leave the details to the reader.
\end{enumerate}
\vspace{0.1in}
We have shown that all lines of $Q(4,q)$ meet $\cM$ in exactly $\frac{q-1}{2}$ points. Therefore $\cM$ is a $\frac{(q-1)}{2}$-ovoid of $Q(4,q)$. The proof is complete.
\qed

\begin{remark}

1. Recall that  the hyperplane $z=0$ meets $Q(4,q)$ in an elliptic quadric $Q^{-}(3,q)$. Using the definition of ${\mathcal M}$, we see that $|{\mathcal M}\cap Q^{-}(3,q)|=|O(1,1,\mu,0)|=\frac{q^2-1}{2}$, which is congruent to $\frac{q-1}{2}\pmod p$. This agrees with the result of Ball \cite{ball} which asserts that an elliptic quadric is incident with $m$ modulo $p$ points of an $m$-ovoid of $Q(4, q)$, where $q$ is a power of $p$. We further conjecture that an elliptic quadric is incident with $m$ modulo $q$ points of an $m$-ovoid of $Q(4, q)$.

2. Using a computer we checked that the subgroup of $\PGO(5,q)$ stabilizing $\cM$ is exactly $A$ when $q=7$ or 11.

\end{remark}

\section{A family of $\frac{(q+1)}{2}$-ovoids of $Q(4,q)$}
%\subsection{Model of $Q(4,q)$ and a prescribed automorphism}
We start this section by commenting that there exist $\frac{(q+1)}{2}$-ovoids of $Q(4,q)$ for all odd $q$.  Cossidente and Penttila \cite{CP} constructed a family of $\frac{(q+1)}{2}$-ovoids ${\mathcal O}$ of the elliptic quadric $Q^{-}(5,q)$ for all odd $q$. Let $H$ be a non-tangent hyperplane of $Q^-(5,q)$. That is, $H\cap Q^-(5,q)$ is a non-singular parabolic quadric in $H$. Then $H\cap {\mathcal O}$ is a $\frac{(q+1)}{2}$-ovoid of $Q(4,q)=H\cap Q^-(5,q)$. In this section, we give a direct construction of $\frac{(q+1)}{2}$-ovoid of $Q(4,q)$ for $q\equiv 1\pmod 4$ by using the same method as in Section 2.

Let $q\equiv 1\pmod{4}$ be a prime power. Let $\gamma$ be a primitive element of $\F_{q^2}$. Furthermore, for positive integers $d\,|\,(q^2-1)$ and  $e\,|\,(q-1)$, let $C_0^{(d,q^2)}$ denote the subgroup of index $d$ of $\F_{q^2}^*$ and put $C_0^{(e,q)}=C_0^{(e(q+1),q^2)}$. Define  $C_{i}^{(d,q^2)}:=\gamma^{i} C_0^{(d,q^2)}$, $0\leq i\leq d-1$, and $C_{i}^{(e,q)}:=\gamma^{i(q+1)} C_0^{(e,q)}$, $0\leq i\leq e-1$.

We will use the following model of  $Q(4,q)$ in this section: Let $V=\F_q\times\F_{q^2}\times\F_{q^2}$, which is viewed as a $5$-dimensional vector space over $\F_q$. We equip $V$ with the following quadratic form
\[
f(x,y,z)=x^2+\tr_{q^2/q}(yz), \; \forall (x,y,z)\in V,
\]
where $\tr_{q^2/q}$ is the trace map from $\F_{q^2}$ to $\F_q$. It is easy to see that $f$ is a non-degenerate quadratic form on $V$, hence the set of zeros of $f$,  $\{(x,y,z)\in V\mid  f(x,y,z)=0\}$, defines a parabolic quadric of $\PG(4,q)$. This is the model that we will be using for $Q(4,q)$ in this section. Formally, we define 
$$Q(4,q)=\{\langle (x,y,z)\rangle\mid 0\neq (x,y,z)\in V, \; f(x,y,z)=0\}.$$
In the rest of this section, in order to simplify notation, we often simply write $(x,y,z)$ for the projective point $\langle(x,y,z)\rangle$ of $\PG(4,q)$. For future use, we also note that the polar form $B$ of $f$ is defined by
\[
B((x_1,y_1,z_1),(x_2,y_2,z_2))= 2x_1x_2+y_1z_2+y_2z_1+y_1z_2^q+y_1^qz_2, \quad \forall (x_1,y_1,z_1), (x_2,y_2,z_2)\in V. 
\]

We now begin to state our construction of $(q+1)/2$-ovoids of $Q(4,q)$. Let $H$ be the cyclic subgroup of order $\frac{q^2-1}{2}$ of $\PGO(5,q)$ consisting of the following elements $T_u$, where $u\in C_0^{(2,q^2)}$. Here $T_u$ is defined by
\[
T_{u}(x,y,z)=(x,yu,zu^{-1}),\quad \forall (x,y,z)\in \PG(V).
\]
Furthermore define two elements $\sigma$ and $\tau$ of order $2$ of $\PGO(5,q)$ as follows:
\begin{align*}
\sigma:\,&(x,y,z)\mapsto (x,y^q,z^q),\quad \forall (x,y,z)\in \PG(V),\\
\tau:  \,&(x,y,z)\mapsto (x,z,y), \quad \forall (x,y,z)\in \PG(V).
\end{align*}
Let $A=\langle H, \sigma, \tau\rangle$ be the subgroup of $\PGO(5,q)$ generated by $H, \sigma$ and $\tau$. We see that $A$ is a subgroup of order $2(q^2-1)$, and it is isomorphic to $C_{\frac{q^2-1}{2}}\rtimes(C_2\times C_2)$.

Let $\delta=\gamma^\frac{q+1}{2}$, where $\gamma$ is a fixed primitive element of $\F_{q^2}$. Then we have $\delta+\delta^q=0$. The group $A$ acts on the set of points of $Q(4,q)$.  Again for a point $P\in\PG(V)$, we use $O(P)$ to denote the orbit in which $P$ lies in. There are four orbits which have a representative having at least one zero coordinate. We list these four orbits below.
\begin{enumerate}
\item $O(0,0,1)$, $O(0,0,\delta)$, each of length $q+1$;
\item $O(0,1,\delta)$ and $O(0,1,\delta^{-1})$, each of length $\frac{q^2-1}{2}$.
\end{enumerate} 

We now give our construction of the $(q+1)/2$-ovoid $\cM$. We define the set $\cM$ to be the union of the following parts:
\begin{enumerate}
\item The orbit $O(0,0,1)$ of length $q+1$;
\item The orbit $O(1,1,-1/2)$ of length $\frac{q^2-1}{2}$;
\item The orbit $O(0,1,\delta)$ of length $\frac{q^2-1}{2}$; 
%\item The elements $\{(1,y,y^{-1}z):\,z\in S,y \in \F_{q^2}^\ast\}\subseteq Q(4,q)$, where 
%\[
%S:=\Big\{z\in \F_{q^\ast}\,|\,z^{q+1}-\frac{1}{4}\in C_{\epsilon}^{(4,q)},z+z^q+1=0\Big\}. 
%\]
\item The set $S=\{(1,y,y^{-1}(w-\frac{1}{2}))\mid w\in C_{\frac{q+1}{2}}^{(2(q+1),q^2)}, \;y \in \F_{q^2}^\ast\}$. 
\end{enumerate}
%By the assumption $(1,y,y^{-1}z)\in Q(4,q)$, $z$ satisfies that $z+z^q+1=0$, 
%i.e., $z=z'-\frac{1}{2}$ for $z'\in C_{\frac{q+1}{2}}^{(q+1,q^2)}$. Then, the condition 
%$z^{q+1}-1/4\in C_{\epsilon}^{(4,q)}$ is reformulated as $-{z'}^2\in C_{\epsilon}^{(4,q)}$. 
Clearly we have
%$|S|=(q-1)/2$ and 
$|{\mathcal M}|=(q^2+1)(q+1)/2$. A couple of remarks are in order. First, the condition $w\in C_{\frac{q+1}{2}}^{(2(q+1),q^2)}$ is equivalent to 
\begin{equation}\label{equiv}
z^{q+1}-\frac{1}{4}\in C_{\epsilon}^{(4,q)} \;{\rm and}\; z^q+z+1=0,\; {\rm where}\; z=w-\frac{1}{2}, 
\end{equation}
where $\epsilon\equiv 1$ or $3$ according as $q\equiv 1$ or $5\pmod 8$. This can be seen as follows. Assume that \eqref{equiv} holds. Then $w^q+w=(z+\frac{1}{2})^q+(z+\frac{1}{2})=z^q+z+1=0$. So $w\in C_{\frac{(q+1)}{2}}^{(q+1,q^2)}$. Furthermore, $z^{q+1}-\frac{1}{4}=(w-\frac{1}{2})^{q+1}-\frac{1}{4}=w^{q+1}$. So by assumption, $w^{q+1}\in C_{\epsilon}^{(4,q)}$, which implies that $w\in C_{\frac{(q+1)}{2}}^{(2(q+1),q^2)}$. The converse is almost the same. We omit the details. With this observation, we comment that $S$ defined above is invariant under the action of $A$.  The set $S$ is clearly invariant under $H$ and $\tau$. To see that $S$ is also invariant under $\sigma$, let $(1,y,y^{-1}(w-\frac{1}{2}))\in S$. Then $\sigma(1,y,y^{-1}(w-\frac{1}{2}))=(1,y^q,y^{-q}(w-\frac{1}{2})^q)=(1,y^q,y^{-q}z^q)$, where $z=w-\frac{1}{2}$. Now $z^q=-z-1=-w-\frac{1}{2}$, and $-w\in C_{\frac{(q+1)}{2}}^{(2(q+1),q^2)}$. We see that $\sigma(1,y,y^{-1}(w-\frac{1}{2}))\in S$. So $S$ is invariant under the action of $A$, hence it is a union of orbits of the action of $A$.

\begin{thm}\label{main2}
The above set $\cM$ is a $\frac{(q+1)}{2}$-ovoid of $Q(4,q)$.
\end{thm}

The proof of the theorem can be done by checking that every line of $Q(4,q)$ meets $\cM$ in exactly $\frac{(q+1)}{2}$ points. The detailed argument is quite similar to that in the proof of Theorem~\ref{main1}. We omit the details.

\begin{remark}
By using a computer we checked that the subgroup of $\PGO(5,q)$ stabilizing the $\frac{(q+1)}{2}$-ovoid $\cM$ above is exactly $A$ when $q=5$ or $9$. In particular, we checked that for $q=5$, there are exactly two inequivalent $\frac{(q+1)}{2}$-ovoids with full automorphism groups of order $48$ in $\PGO(5,q)$, and there are two inequivalent $\frac{(q+1)}{2}$-ovoids from the hyperplane sections of the Cossidente-Penttila ovoid with full automorphism groups of order $48$. Hence our construction in this case does not produce new $\frac{(q+1)}{2}$-ovoids. For larger $q$, the amount of computations involved in checking inequivalence becomes too large to handle.
\end{remark}

\section{Remarks and Open Problems}
For $q$ even, the generalized quadrangle $Q(4,q)$ contains $m$-ovoids for every value of $m$ (cf. \cite{OW3q}). For $q$ odd, $Q(4, q)$ contains $m$-ovoids for $m=1$, $m = (q + 1)/2$ and $m = q$.  It is not known in general for what other values of $m$ the generalized quadrangle $Q(4,q)$ contains $m$-ovoids. In this paper, by construction we have shown that $Q(4,q)$ contains $\frac{(q-1)}{2}$-ovoids when $q\equiv 3\pmod 4$. For small prime powers $q\equiv 3\pmod 4$,  we have checked that the stabilizer in $\PGO(5,q)$ of the newly constructed $\frac{(q-1)}{2}$-ovoid has size $2(q^2-1)$ (that is, the full automorphism group in $\PGO(5,q)$ of the newly constructed $\frac{(q-1)}{2}$-ovoid is exactly the group $A$ that we prescribed). As for future research, it would be interesting to generalize the examples of $\frac{q-1}{2}$-ovoid of $Q(4,q)$ (with $q=5$ or $9$) in Example 5 of \cite{bambergCombin} into an infinite family. We could not see any general pattern for the prescribed automorphism groups in those examples. Secondly, Ball \cite{ball} proved that an $m$-ovoid of $Q(4,q)$ meets every three-dimensional elliptic quadric in $m$ modulo $p$ points, where $q$ is a power of $p$. We conjecture that an $m$-ovoid of $Q(4,q)$ meets every three-dimensional elliptic quadric in $m$ modulo $q$ points. There are some results related to this conjecture in \cite{ballstorme} and \cite{chandler}.


\begin{thebibliography}{99}

\bibitem{ball} S. Ball, On $m$-ovoids of parabolic quadrics, unpublished notes.

\bibitem{ballstorme} S. Ball, P. Govaerts, L. Storme, On ovoids of parabolic quadrics, {\it Des. Codes  Cryptogr.} {\bf 38} (2006), 131-145.

\bibitem{bambergjcta} J. Bamberg, S. Kelly, M. Law, T. Penttila, Tight sets and $m$-ovoids of finite polar spaces, {\it J. Combin. Theory} (A) {\bf 114} (2007), 1293--1314.

\bibitem{bambergCombin} J. Bamberg, M. Law, T. Penttila, Tight sets and $m$-ovoids of generalised quadrangles, {\it Combinatorica}, {\bf 29} (2009), 1--17.

\bibitem{bambergDevillers} J. Bamberg, A. Devillers, J. Schillewaert, Weighted intriguing sets of finite generalised quadrangles, {\it J. Alg. Combin.} {\bf 36} (2012), 149--173.

%\bibitem{fans}A. E. Brouwer and H. A. Wilbrink. Ovoids and fans in the generalized quadrangle $Q(4, 2)$, %Geom. Dedicata, {\bf 36} (1990), 121--124.

%\bibitem{}J. Bamberg, J. De Beule, and F. Ihringer, Weighted intriguing sets in finite polar spaces, arxiv 2015.

\bibitem{cameron}P. J. Cameron, Automorphisms of graphs, {\it Topics in algebraic graph theory}, {\bf 102} (2004), 137--155.

\bibitem{chandler} D. B. Chandler, The sizes of the intersections of two unitals in $\PG(2,q^2)$, {\it Finite Fields Appl.}  {\bf 25} (2014), 255-–269.

\bibitem{OW3q} A. Cossidente, C. Culbert, G. L. Ebert,  G. Marino, On $m$-ovoids of $W_3(q)$, {\it Finite Fields Appl.} {\bf 14} (2008), 76--84.

\bibitem{CP} A. Cossidente, T. Penttila, Hemisystems of the Hermitian surfaces, {\it J. London Math. Soc. (2)} {\bf 72 (3)} (2005), 731--741.

\bibitem{jande} J. De Beule,  J. Demeyer, K. Metsch, M. Rodgers, A new family of tight sets in
$\cQ^+(5,q)$, in print {\it Des. Codes Cryptogr.}, DOI 10.1007/s10623-014-0023-9

\bibitem{drudge} K. Drudge, {\it Extremal sets in projective and polar spaces}, PhD theses, The University of Western Ontario, 1998.

\bibitem{er} P. Erd\H{o}s, A. R\'enyi, Asymmetric graphs, {\it Acta Math. Acad. Sci. Hungar.} {\bf 14} (1963), 295--315.

\bibitem{fmx} T. Feng, K. Momihara, Q. Xiang, Cameron-Liebler line classes with parameter $x=\frac{q^2-1}{2}$, {\it J. Combin. Theory} (A) {\bf 133} (2015), 307--338.

%\bibitem{pavese}Francesco Pavese, Geometric constructions of two-character sets, Discrete Mathematics %338 (2015) 202-208

%\bibitem{}Bart De Bruyn, Intriguing Sets of Points of $Q(2n, 2)\setminus Q^+(2n-1,2)$, Graphs and %Combinatorics (2012) 28:791-805

%\bibitem{fieldred}Kelly, Constructions of intriguing sets of polar spaces from field reduction and derivation, %Des Codes Crypt (2007) 43:1-8.

%\bibitem{w5q} Antonio Cossidente, Francesco Pavese, Intriguing sets of $W(5, q)$, $q$ even, J. Combin. Theory Ser. A 127 (2014) 303-313

%\bibitem{bier}J. Bierbrauer, and Y. Edel, A family of 2-weight codes related to BCH-codes, J Combin Des 5 %(1997), 391-396.

%\bibitem{cycl}Gustavo A. Fernandez-Alcober, Rugare Kwashira, Luis Martinez, Cyclotomy over products of f%inite fields and combinatorial applications, European Journal of Combinatorics 31 (2010) 1520-1538

%\bibitem{konert} Konert, Constructing two-weight codes with prescribed groups of
%automorphisms, Discrete Applied Mathematics 155 (2007) 1451 -1457

%\bibitem{GrSRGdata1}Iliya Bouyukliev, Veerle Fack, Wolfgang Willems, Joost Winne,
%Projective two-weight codes with small parameters and their corresponding graphs, Des Codes Crypt (2006) %41:59-78

%\bibitem{}Cossidente, Embeddings of $U_n(q^2)$ and symmetric SRG, J Combin Designs 18: 248-253, 2010

%\bibitem{}Antonio Cossidente, Hendrik Van Maldeghem, The simple exceptional group $G_2(q)$, $q$ even, %and two-character sets, Journal of Combinatorial Theory, Series A 114 (2007) 964-969.

%\bibitem{}A. DeWispelaere1, H. Van Maldeghem, Some new two-character sets in $\PG(5, q^2)$ and a %distance-2 ovoid in the generalized hexagon $H(4)$, Discrete Mathematics 308 (2008) 2976-2983

%\bibitem{}Cossidente, King, Some two-character sets, Des. Codes Cryptogr. (2010) 56:105, C113.

%\bibitem{}The geometry of some two-character sets , Des. Codes Cryptogr. (2008) 46:231-241.

%\bibitem{CP}Cossidente, Penttila, Two-Character Sets Arising from Gluings of Orbits, Graphs and %Combinatorics (2013) 29:399-406.

%\bibitem{GC15}Francesco Pavese, Geometric constructions of two-character sets, , Discrete Mathematics %338 (2015) 202-208.

\bibitem{LN97}
R. Lidl, H. Niederreiter, {\it Finite Fields}, Cambridge
Univ. Press, 1997.
\bibitem{payne} S. E. Payne, Tight pointsets in finite generalized quadrangles, in: Eighteenth Southeastern International Conference on Combinatorics, Graph Theory, and Computing, Boca Raton, FL, 1987, {\it Congr. Numer.}  {\bf 60} (1987), 243--260.

\bibitem{morgan} M. Rodgers,  {\it On some new examples of Cameron-Liebler line classes}, PhD thesis, University of Colorado, 2012.

\bibitem{ST1} E. E. Shult, J. A. Thas, m-systems of polar space, {\it J. Combin. Theory} (A), {\bf 68} (1994), 184--204.

%\bibitem{ST2} E. E. Shult, J.A. Thas, Constructions of polygons from buildings, PLMS (3), 71(2): 397-440, %1995.

%\bibitem{HQ} Hamilton, Quinn, ``m-systems of polar spaces and maximal arcs in projective planes." BBMS, %2000.

%\bibitem{HM} N. Hamilton and R. Mathon. Existence and non-existence of m-systems
%of polar spaces. EJC  2001.

%\bibitem{BKM}Jan De Beule, A. Klein, K. Metsch, Substructures of finite classical polar spaces, in: Current %research topics in Galois geometries (eds: %De Beule, L. Storme).

%\bibitem{HT}Hirschfeld, Thas, Open problems in finite projective spaces, Finite Fields and Their Applications, 32(2015), 44-81.

%\bibitem{W3qOvoid}J. A. Thas, On $4$-gonal configurations, Geom. Dedicata, 2: 317-326, 1973.

\bibitem{thaslaa} J. A. Thas, Interesting pointsets in generalized quadrangles and partial geometries, {\it Linear Algebra Appl.} {\bf 114/115} (1989), 103--131.

\bibitem{wormald} N. C. Wormald, Models of random regular graphs, In: Surveys in Combinatorics, 1999 (ed. J. D. Lamb and D. A. Preece), {\it London Math. Soc. Lecture Notes Series} {\bf 267}, Cambridge University Press, Cambridge, 1999, 239-298. 

%\bibitem{packing}Penttila, Williams, Regular Packings of $PG(3,q)$, Europ. J. Combinatorics (1998) 19, %713-720

\end{thebibliography}
\end{document}